*LA MARAVILLOSA FUNCION Y ECUACION CUADRATICA*

CAMPO ELÍAS GONZALEZ PINEDA.
cegp@utp.edu.co

*Introducción.*

En este artículo presentamos un estudio general de la función cuadrática. Veremos lo importante de esta función y ecuación en todo el ámbito de las ciencias.

Esperamos mostrar al lector la importancia de ciertos contenidos matemáticos que a veces pasan desapercibidos por lo simple y su grado de complejidad. Sin embargo, veremos como algo que parece tan simple en realidad no lo es. El lector podrá encontrar muchas otras aplicaciones de la función cuadrática. Solo queremos resaltar la importancia de esta función.

*Contenido.*

De todas las funciones polinómicas no lineales, quizá la más importante es la función cuadrática dada por:

$$f(x) = x^2 \quad (1)$$

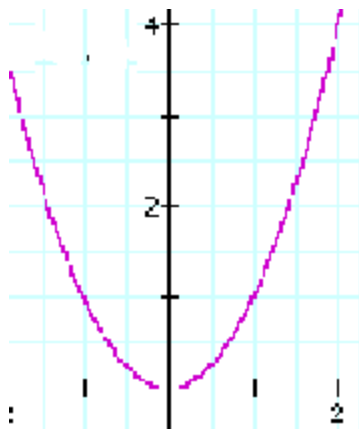

La ecuación la podemos multiplicar por $a$ y sumarle un $b$. Es decir (1) se puede escribir en la forma:

$$f(x) = ax^2 + b, a \neq 0$$

Cuya gráfica abre hacia arriba si $a > 0$ y hacia abajo si $a < 0$. Además, la gráfica está traslada $b$ unidades. El caso más general viene dado por

$$f(x) = ax^2 + bx + c$$

Que se reduce al caso anterior, completando cuadrado

$$f(x) = a\left(x + \frac{b}{2a}\right)^2 + c - \frac{b^2}{4a}$$

Vemos que el vértice tiene coordenadas

$$\left(-\frac{b}{2a}, \frac{4ac - b^2}{4a}\right)$$

La función cuadrática tiene muchas aplicaciones, por ejemplo lanzamiento de proyectiles, crecimiento de poblaciones, ingresos, etc. Por ahora centramos la atención en la ecuación cuadrática

$$ax^2 + bx + c = 0, a \neq 0 \ .$$

Sin pérdida de generalidad podemos considerar la ecuación

$$x^2 + bx + c = 0$$

que por la fórmula cuadrática tiene como solución a

$$x = \frac{-b \pm \sqrt{b^2 - 4c}}{2}$$

Como sabemos las raíces de la ecuación cuadrática pueden ser: reales e iguales, reales distintas, o complejas conjugadas.

La pregunta que surge es ¿qué importancia tiene esta ecuación (y función) en las Matemáticas y su aplicación al mundo real? Esta pregunta la contestaremos en parte de aquí en adelante. Antes mostramos unas particularidades de los ceros de la ecuación cuadrática.

### Curiosidades de la ecuación cuadrática.

Sean $p > 0, q > 0$ reales positivos. Estos números generan cuatro ecuaciones cuadráticas especiales a saber

$$\begin{array}{ll} x^2 + px + q = 0 \ (a), & x^2 - px + q = 0 \ (b) \\ x^2 + px - q = 0 \ (c), & x^2 - px - q = 0 \ (d) \end{array}$$

Nótese que si $r_1, r_2$ son las raíces de (a) entonces $-r_1, -r_2$ son las raíces de (b) y viceversa. Esto es lógico porque una ecuación se obtiene de la otra cambiando $x$ por $-x$. Un resultado similar se tiene para las ecuaciones (c) y (d). Por la formula cuadrática encontramos para el caso (b) y (d) que las raíces vienen dadas por:

$$x_b = \frac{p \pm \sqrt{p^2 - 4q}}{2}$$

$$x_d = \frac{p \pm \sqrt{p^2 + 4q}}{2}$$

Nótese que para el primer caso las raíces son reales y distintas si $p^2 > 4q$ y son complejas conjugadas si $p^2 < 4q$. En ambos casos si $p^2 - 4q$ no es un cuadrado perfecto, las raíces están en el cuerpo cuadrático $Q(\sqrt{p^2 - 4q})$. Para el segundo caso las raíces siempre son reales distintas, y están en el cuerpo cuadrático $Q(\sqrt{p^2 + 4q})$ en el caso de que la cantidad subradical no sea un cuadrado perfecto.

De otro lado notemos que

$$x_b + k = \frac{p + 2k \pm \sqrt{p^2 - 4q}}{2}$$

Números que son las raíces de la ecuación

$$x^2 - (p + 2k)x + k^2 + pk + q = 0$$

Es interesante notar que aparece la otra ecuación y que además $(k^2 + pk + q)' = p + 2k$.

De manera similar

$$x_a + k = \frac{-p + 2k \pm \sqrt{p^2 - 4q}}{2}$$

Son las raíces de la ecuación

$$x^2 - (-p + 2k)x + k^2 - pk + q = 0$$

Esto se obtiene cambiando $p$ por $-p$.
De manera similar

$$x_d + k = \frac{p + 2k \pm \sqrt{p^2 + 4q}}{2}$$

Que son las raíces de la ecuación

$$x^2 - (p + 2k)x + k^2 + pk - q = 0$$

Similar al caso anterior aparece la otra ecuación y que además $(k^2 + pk - q)' = p + 2k$.
Ahora, cambiando $p$ por $-p$ obtenemos la ecuación

$$x^2 - (-p + 2k)x + k^2 - pk - q = 0$$

Haciendo $p = q$ tenemos los siguientes casos particulares.
$$x^2 - px - p = 0$$

Encontramos
$$r = \frac{p \pm \sqrt{p^2 + 4p}}{2},$$

Ahora si sumamos 1 a estas raíces tenemos
$$r + 1 = \frac{p + 2 \pm \sqrt{p^2 + 4p}}{2},$$

que son los ceros de la ecuación
$$x^2 - (p + 2)x + 1 = 0$$

De manera similar para la ecuación
$$x^2 + px - p = 0,$$
Vemos que las raíces son
$$r = \frac{-p \pm \sqrt{p^2 + 4p}}{2},$$
Luego
$$r + 1 = \frac{-p + 2 \pm \sqrt{p^2 + 4p}}{2}$$
Que son las raíces de la ecuación
$$x^2 - (-p + 2)x - (4p - 1) = 0$$

Para la ecuación
$$x^2 - px + p = 0,$$
Tenemos
$$r = \frac{p \pm \sqrt{p^2 - 4p}}{2}$$
Sumando uno tenemos

$$r + 1 = \frac{p + 2 \pm \sqrt{p^2 - 4p}}{2}$$

Que son los ceros de la ecuación

$$x^2 - (p+2)x + 4p + 1 = 0$$

Para la ecuación

$$x^2 + px + p = 0,$$

Encontramos

$$r = \frac{-p \pm \sqrt{p^2 - 4p}}{2} \leftrightarrow r + 1 = \frac{-p + 2 \pm \sqrt{p^2 - 4p}}{2}$$

Que son las raíces de la ecuación

$$x^2 - (-p + 2)x + 1 = 0$$

Nótese la relación entre las raíces. Además, restando $p$ a las raíces de una se obtienen las de la otra.

## *Cuerpos cuadráticos*

Sea $K$ un cuerpo y supongamos que la ecuación

$$x^2 = m \quad (2)$$

no tiene solución en $K$. Construyamos el conjunto

$$K(\sqrt{m}) = \{a + b\sqrt{m} : a, b \in K\}$$

Si $z = a + b\sqrt{m}$, $w = c + d\sqrt{m}$ son elementos de $K(\sqrt{m})$ definimos,

1. $z = w$ si y solo si $a = c, b = d$.
2. $z + w = (a + c) + (b + d)\sqrt{m}$
3. $z \cdot w = (ac + bdm) + (bc + ad)\sqrt{m}$

   No es difícil mostrar que $(K(\sqrt{m}), +, \cdot)$ tiene estructura de cuerpo (llamado cuerpo cuadrático generado por m asociado a $K$) y además:

1. La ecuación (2) tiene solución en $K(\sqrt{m})$.
2. Si $W$ es un cuerpo en el que (2) tiene solución entonces,
$$K(\sqrt{m}) \subseteq W$$
3. $K \subseteq K(\sqrt{m})$
4. $K(\sqrt{m}) \equiv K \times K$ ($\equiv$ *es isomorfo.*)

Es bueno anotar que en $K(\sqrt{m})$ la solución de la ecuación (2) es $x = \pm\sqrt{m}$.

De otra parte si hacemos $j = \sqrt{m}$ tenemos que $z = a + bj$. Al número $a$ se le llama parte real de $z$ y al número $b$ parte imaginada de $z$. El conjugado de $z$ se denota y define $\bar{z} = a - bj$. Pero porqué hemos centrado la razón en los cuerpos cuadráticos? pues precisamente porque son generados por una ecuación cuadrática. Ahora, sabemos que en los números reales ($R$) la ecuación $x^2 = -1$ no tiene solución, por lo que tenemos el cuerpo cuadrático $C = R(\sqrt{-1})$ llamado el cuerpo de los números complejos. Además este resultado no da mucha información, por ejemplo, entre los números reales y los números complejos no hay cuerpos, mientras que entre los números racionales y los reales hay infinitos cuerpos. Existen infinitos cuerpos entre los racionales y los complejos que no están contenidos en los reales.

## Ecuación diferencial ordinaria de orden dos con coeficientes constantes.

Del curso básico de ecuaciones diferenciales ordinarias se conoce la ecuación con coeficientes constantes

$$ay'' + by' + cy = 0 \quad a \neq 0 \quad (3)$$

Buscando soluciones de la forma $y = e^{mx}$ se encuentra que $m$ es un cero de la ecuación

$$am^2 + bm + c = 0$$

Llamada ecuación característica asociada a (3).

Sabemos que la ecuación (3) modela muchos fenómenos físicos: movimiento armónico simple, amortiguado, circuitos en serie, el péndulo etc. Solo por mencionar algunos. Una consecuencia importante, es la relación en este caso de la ecuación de segundo grado con la función exponencial. Más adelante mencionaremos esto. Recuérdese que la ecuación

$$ax^2 y'' + bxy' + cy = 0$$

Puede escribirse como una ecuación de coeficientes constantes de segundo orden.

## Método de Frobenius.

Cuando se aplica el método de Frobenius a la ecuación $y'' + P(x)y + Q(x)y = 0$ la forma de las soluciones depende de una ecuación cuadrática de la forma $r^2 + br + c = 0$ llamada ecuación indicial.

## *La ecuación cuadrática y su derivada.*

Recordemos que si $ax^2 + bx + c = 0$ entonces por la formula cuadrática,

$$x = \frac{-b \pm \sqrt{b^2 - 4ac}}{2a} = \frac{-b \pm h}{2a}$$

Sea $x_1 = \frac{-b+h}{2a}$ entonces, $h = 2ax_1 + b = f'(x_1)$. Es decir, la raíz cuadrada del discriminante de la ecuación, es la derivada de la función cuadrática evaluada en el cero dado. Es decir, $\Delta = \left(f'(x_1)\right)^2$. Además nótese que el segundo cero es $x_2 = -\frac{ax_1+b}{a}$

## *La ecuación cuadrática y la sucesión de Fibonacci y Teoría de grupos.*

Consideremos los siguientes casos particulares de la ecuación cuadrática:

### *Caso I: Ecuación*

$$x^2 = x + 1 \quad (1)$$

Los ceros de esta ecuación son el número áureo y su conjugado. $r = \frac{1 \pm \sqrt{5}}{2}$. Es de anotar que el número áureo aparece de manera natural en el universo. La ecuación anterior la podemos escribir en la forma

$$x = 1 + \frac{1}{x}$$

Es curioso que la expresión $\frac{1}{x}$ sea la parte decimal de $x$.

Multiplicando por $x^{n-1}$ encontramos

$$x^n = x^{n-1} + x^{n-2} \leftrightarrow x^n - x^{n-1} = x^{n-2}$$

Sumando desde 1 hasta n obtenemos

$$\sum_{k=1}^{n}(x^k - x^{k-1}) = \sum_{k=1}^{n} x^{k-2}$$

Esto es $\sum_{k=1}^{n} x^k = x^{n+2} - x^2$

De otro lado la sucesión.

$$x^n = x^{n-1} + x^{n-2}$$

es una generalización de la sucesión de Fibonacci. No es difícil mostrar en la ecuación anterior(o de la ecuación (1)) que

$$x^n = F_n x + F_{n-1}$$

donde $F_n$ es el enésimo número de Fibonacci. Enlistados los coeficientes quedan

| $F_n$ | $F_{n-1}$ |
|-------|-----------|
| 1     | 1         |
| 2     | 1         |
| 3     | 2         |
| 5     | 3         |
| 8     | 5         |
| 13    | 8         |
| ⋮     | ⋮         |

Vemos en este caso como la ecuación cuadrática está relacionada con la sucesión de Fibonacci y además sus raíces generan a su vez una sucesión generalizada de Fibonacci. Recordemos que la sucesión de Fibonacci se define por la fórmula de recurrencia:

$$f_n = f_{n-1} + f_{n-2}, \quad f_0 = 1, f_1 = 1$$

Sumando encontramos

$$\sum_{k=1}^{n} f_{k-2} = \sum_{k=1}^{n}(f_k - f_{k-1}) = f_n - f_0 = f_n - 1$$

Observemos que para la ecuación (1) $x = \frac{1 \pm \sqrt{5}}{2}$

Consideremos el caso general $x^b = x + 1$. Haciendo

$$x^{b+k-1} = x^k + x^{k-1}$$

No es difícil ver que

$$\sum_{k=1}^{n} x^k = \frac{x^b - 1}{x - 1}(x^n - 1) = \frac{x}{x - 1}(x^n - 1)$$

La prueba es por inducción, claramente para $n = 1$ el resultado se tiene. Supongámoslo para $n$ y veamos que se cumple para $n + 1$. En efecto,

$$\sum_{k=1}^{n+1} x^k = \sum_{k=1}^{n} x^k + x^{n+1}$$

$$= \frac{x}{x-1}(x^n - 1) + x^{n+1} = \frac{x^{n+1} - x + x^{n+2} - x^{n+1}}{x-1} = \frac{x}{x-1}(x^{n+1} - 1)$$

Nótese que simplemente el resultado es consecuencia de la suma geométrica.

*Caso II: Ecuación*

$$x^2 = -x + 1 \quad (2)$$

Multiplicando por $x^{k-2}$ obtenemos,

$$x^k = x^{k-2} - x^{k-1} \quad (**)$$

Sumando tenemos,

$$\sum_{k=1}^{n} x^k = \frac{1}{x} - x^{n-1} = \frac{1 - x^n}{x} = (x+1)(1 - x^n)$$

De la ecuación (**) o de (2) encontramos

$$x^n = (-1)^{n+1} F_n x + (-1)^n F_{n-1}$$

En este caso la solución de la ecuación es $x = \frac{-1 \pm \sqrt{5}}{2}$

Resolvamos el ahora

$$x^b = 1 - x \leftrightarrow x = x^{1-b}(1 - x)$$

De donde

$$x^k = x^{1-b}(x^{k-1} - x^k)$$

Sumando

$$\sum_{k=1}^{n} x^k = x^{1-b}(1 - x^n) = \frac{x(1 - x^n)}{x^b}$$

***Caso III Ecuación***

$$x^2 = x - 1 \quad (3)$$

De tres encontramos que

$$x^{2+6t} = -x^{5+6t} = x - 1$$
$$x^{4+6t} = -x^{7+6t} = -x$$
$$x^{6+6t} = -x^{3+6t} = 1$$

Para $t \geq 0$. Es fácil ver que el conjunto
$$T = \{x^2, x^3, x^4, x^5, x^6, x^7\} = \{x^2, x^3, x^4, -x^2, -x^3, -x^5\}$$

es un grupo abeliano con el producto de orden 6. Ahora,

$$x^2 = x - 1 \leftrightarrow x = 1 - x^{-1}$$

Si multiplicamos por $x^{k-1}$, es decir, $x^k = x^{k-1} - x^{k-2}$ entonces

$$\sum_{k=1}^{n} x^k = x^{n-1} - x^{-1} = \frac{x^n - 1}{x}$$

Tenemos que
1. Si $k = 2 + 6t$ (o $k = 5 + 6t$)
$$\sum_{k=1}^{n} x^k = \frac{x - 2}{x} \quad \left(o = -\frac{x + 1}{x}\right)$$
2. Si $k = 4 + 6t$ (o $k = 7 + 6t$)
$$\sum_{k=1}^{n} x^k = \frac{-x - 1}{x} \quad \left(o = \frac{x - 1}{x} = x\right)$$
3. Si $k = 6 + 6t$ (o $k = 3 + 6t$)

$$\sum_{k=1}^{n} x^k = 0 \quad \left(o = \frac{-2}{x}\right)$$

Nótese que en este caso

$$x = \frac{1 \pm \sqrt{3}i}{2}$$

En términos complejos:

$$x = e^{\frac{\pi}{3}i} \quad o \quad e^{\frac{-\pi}{3}i}$$

Considere el caso general

$$x^b = x - 1$$

Así,

$$x^k = x^{1-b}(x^k - x^{k-1})$$

Sumando

$$\sum_{k=1}^{n} x^k = x^{1-b} \sum_{k=1}^{n} (x^k - x^{k-1}) = x^{1-b}(x^n - 1) \quad (*)$$

El lado derecho de la ecuación (*), se puede escribir de distintas maneras según reemplazar $x^{-b}$, algunas de ellas son:

- $x^{-b} = \frac{1}{x-1}$, tenemos
$$\sum_{k=1}^{n} x^k = \frac{x(x^n - 1)}{x - 1}$$
- $x^{-b} = x^{1-b} - x$ luego
$$\sum_{k=1}^{n} x^k = x^{1-b}(x^n - 1) = x^{n+1-b} - x^{1-b} = x^{n+1-b} - x^{-b} - x$$
- $x^{-b} = \frac{1}{x-1}$ luego
$$\sum_{k=1}^{n} x^k = x^{1-b}(x^n - 1) = x^{n+1-b} - x^{1-b} = x^{n-b} - \frac{x}{x - 1}$$

### *Caso IV. Ecuación*

$$x^2 = -x - 1 \quad (4)$$

Multiplicando la ecuación (4) por $x$ e iterando el proceso encontramos que
$$x^2 = -x - 1, x^3 = 1, x^4 = x, x^5 = x^2 \ldots.$$

Y así el conjunto
$$T = \{x^2, x^3, x^4\}$$
tiene estructura de grupo abeliano y es de orden 3. En general se tiene que
$$x^{2+3t} = x^2, \qquad x^{3+3t} = x^3 = 1, \ x^{4+3t} = x$$
De la fórmula
$$\sum_{k=1}^{n} x^k = \frac{x(x^n - 1)}{x - 1}$$
Encontramos lo siguiente:
1. $n = 2 + 3t, \ x^n = x^2$ luego
$$\sum_{k=1}^{n} x^k = \frac{x(x^2 - 1)}{x - 1} = \frac{x^3 - x}{x - 1} = -1.$$
2. $n = 3t, \ x^n = x^3$ luego
$$\sum_{k=1}^{n} x^k = \frac{x(x^3 - 1)}{x - 1} = \frac{x(1 - 1)}{x - 1} = 0.$$
3. $n = 4 + 3t, \ x^n = x$ luego
$$\sum_{k=1}^{n} x^k = \frac{x(x - 1)}{x - 1} = x.$$

En caso $x = \frac{-1 \pm \sqrt{3}i}{2}$, es decir, $x_1 = e^{\frac{2\pi i}{3}}, \ x_2 = e^{-\frac{2\pi i}{3}}$

***Relación entre el número áureo y pi y la ecuación cuadrática.***

La relación entre $\varphi$ y $\pi$ se deduce de la identidad trigonométrica

$$sen(5\theta) = 5sen(\theta) - 20sen^3(\theta) + 16sen^5(\theta)$$
Haciendo $\theta = \frac{\pi}{5}$. Resulta una ecuación cuadrática, de donde se obtiene
$$\cos^2\left(\frac{\pi}{5}\right) = \frac{\varphi^2}{4}, sen^2\left(\frac{\pi}{5}\right) = \frac{\overline{\varphi^2}}{4} + \frac{1}{4}$$

Sacando raíz cuadrada se obtiene el coseno y el seno (raíz positiva.) Nótese que
$$\frac{\varphi^2}{4}+\frac{\overline{\varphi^2}}{4}+\frac{1}{4}=1 \leftrightarrow \varphi^2+\overline{\varphi^2}=3$$
El último resultado es de especial interés y lo bautizamos.

**Definición:** La ecuación de la creación viene dada por
$$\varphi^2+\overline{\varphi^2}=3.$$

El resultado anterior se puede generalizar, recuérdese que $\phi$ es cero de la ecuación cuadrática
$$x^2-x-1=0$$
cuya solución está en el cuerpo $Q(\sqrt{5})$.

**Generalización:** De igual forma podemos considerar los números
$$\Phi=\frac{1+\sqrt{m}}{2}, \quad \overline{\Phi}=\frac{1-\sqrt{m}}{2}$$
Que son soluciones de la ecuación
$$x^2-x-\frac{m-1}{4}=0$$
válidas para todo $m$, pero nos interesa el caso en el que el número $m$ sea un entero no cuadrado perfecto. Nótese que
$$\Phi^2+\overline{\Phi^2}=\frac{m+1}{2}$$

Si hacemos
$$\cos(\theta)=\frac{\Phi}{2}$$
Encontramos que
$$sen^2(\theta)=1-\left(\frac{\Phi}{2}\right)^2=\frac{15-m+\left(1-\sqrt{m}\right)^2-m-1}{16}=\frac{7-m}{8}+\frac{\overline{\Phi^2}}{4}$$
Es decir,
$$\frac{\Phi^2}{4}+\frac{7-m}{8}+\frac{\overline{\Phi^2}}{4}=1$$
Los valores permitidos para $m$ son 0, 1, 2, 3, 4, 5, 6, 7, 8,9.
Por último notemos que
$$2\cos\left(\frac{\pi}{5}\right)=\frac{1+\sqrt{5}}{2}, \quad 2\cos\left(\frac{2\pi}{5}\right)=\frac{-1+\sqrt{5}}{2}$$

***Un caso especial.***

Si hacemos
$$x_1 = \frac{1 + 2k + \sqrt{m}}{2}, \quad x_2 = \frac{1 + 2k - \sqrt{m}}{2}$$
Encontramos
$$x_1 + x_2 = 2k + 1$$
$$x_1 \cdot x_2 = k^2 + k + \frac{1-m}{4}$$
y obtenemos la ecuación
$$x^2 - (2k+1)x + k^2 + k + \frac{1-m}{4} = 0 \ (\Delta)$$
De igual forma si hacemos
$$x_1 = \frac{-1 + 2k + \sqrt{m}}{2}, \quad x_2 = \frac{-1 + 2k - \sqrt{m}}{2}$$
Encontramos
$$x_1 + x_2 = 2k - 1$$
$$x_1 \cdot x_2 = k^2 - k + \frac{1-m}{4}$$
Obtenemos la ecuación
$$x^2 - (2k-1)x + k^2 - k + \frac{1-m}{4} = 0 \ (\Delta_1)$$
Las ecuaciones $(\Delta_1)$ y $(\Delta)$ tienen una relación muy fuerte, y su forma es realmente similar.
De manera análoga si suponemos $m < 0$ y hacemos
$$x_1 = \frac{1 + 2k + \sqrt{m}i}{2}, \quad x_2 = \frac{1 + 2k - \sqrt{m}i}{2}$$
Obtenemos la ecuación
$$x^2 - (2k+1)x + k^2 + k + \frac{1+m}{4} = 0$$
Si
$$x_1 = \frac{-1 + 2k + \sqrt{m}i}{2}, \quad x_2 = \frac{-1 + 2k - \sqrt{m}i}{2}$$

$$x^2 - (2k-1)x + k^2 - k + \frac{1+m}{4} = 0$$
Y vemos su relación con la ecuación anterior.

## *Ecuación cuadrática y teoría de grupos.*

Recordemos que las soluciones de la ecuación
$$x^2 = -x - 1 \ (4)$$

Genera el grupo
$$T = \{x^2, x^3, x^4\}$$
Con el producto, además este grupo es cíclico. Recordemos que
$$x^2 = -x - 1, x^3 = 1, x^4 = x, x^5 = x^2 \ldots$$

También la ecuación
$$x^2 = x - 1 \quad (3)$$
genera el grupo
$$T = \{x^2, x^3, x^4, x^5, x^6, x^7\}$$
Recordemos que:
$$x^{2+6t} = -x^{5+6t} = x - 1$$
$$x^{4+6t} = -x^{7+6t} = -x$$
$$x^{6+6t} = -x^{3+6t} = 1$$

Para $t \geq 0$. Algunas propiedades de conmutatividad en grupo quedan definidas por una ecuación cuadrática:
1. Sea $(G, \cdot)$ un grupo con identidad $e$. Si para todo $a \in G$ se tiene que $a^2 = e$ el grupo es abeliano.
2. Sea $(G, \cdot)$ un grupo, se cumple que $(ab)^2 = a^2 b^2$ para todo los elementos de $G$ entonces, $G$ es abeliano.

Pero quizá la relación más importante de la ecuación cuadrática con los grupos aparezca en la teoría de congruencias. La cual a su vez nos dará una relación con los números primos, relación de por sí bien interesante, ya que como veremos es un caso especial del teorema de Pitágoras.

## *Ecuación cuadrática y teoría de congruencias.*

Supongamos que se quiere resolver la ecuación cuadrática
$$ax^2 + bx + c \equiv 0 \bmod p$$
equivale a resolver la ecuación
$$u^2 \equiv r \bmod p \ , u = 2ax + b, r = 4ac - b^2 \qquad (t')$$

con $r, p$ primos relativos.

Recuérdese el resultado importante de residuos cuadráticos:

La ecuación $x^2 \equiv -1 \bmod p$ tiene solución si y solo si $p$ es un primo de la forma $4t + 1$. Este resultado equivale a decir que existen enteros $a, b$ de tal manera que

$$a^2 + b^2 = p$$

Como caso particular:

$$1 + 2^2 = 5$$

$$2^2 + 3^2 = 13$$

$$1 + 4^2 = 17$$

## *La función cuadrática y la geometría*

La función cuadrática aparece con mucha frecuencia en la geometría, en particular en las áreas de las figuras geométricas. Por ejemplo, el área de un cuadrado es su lado al cuadrado, el área de un círculo es $\pi r^2$, etc. De hecho en uno de los orígenes del número de oro cuando se hace la proporción

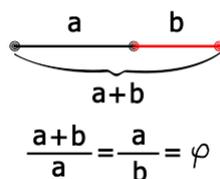

Recuérdese que el área de una esfera de radio $r$ es $A = 4\pi r^2$ y que su volumen es $V = \frac{4}{3}\pi r^3 = \pi r^2 \times \frac{4}{3}r = 4\pi r^2 \times \frac{1}{3}r$ .

| Nombre | Área de una cara | Área total | Apotema | Volumen |
|---|---|---|---|---|
| Tetraedro | $\dfrac{a^2}{4}\cdot\sqrt{3}$ | $a^2\cdot\sqrt{3}$ | $\dfrac{a}{12}\cdot\sqrt{6}$ | $\dfrac{a^3}{12}\cdot\sqrt{2}$ |
| Octaedro | $\dfrac{a^2}{4}\cdot\sqrt{3}$ | $2\cdot a^2\cdot\sqrt{3}$ | $\dfrac{a}{6}\cdot\sqrt{6}$ | $\dfrac{a^3}{3}\cdot\sqrt{2}$ |
| Icosaedro | $\dfrac{a^2}{4}\cdot\sqrt{3}$ | $5\cdot a^2\cdot\sqrt{3}$ | $\dfrac{a}{2}\sqrt{\dfrac{7+3\cdot\sqrt{5}}{6}}$ | $\dfrac{5\cdot a^3}{6}\cdot\sqrt{\dfrac{7+3\cdot\sqrt{5}}{2}}$ |
| Hexaedro | $a^2$ | $6\cdot a^2$ | $\dfrac{a}{2}$ | $a^3$ |
| Dodecaedro | $\dfrac{5}{4}\cdot a^2\cdot\sqrt{\dfrac{5+2\cdot\sqrt{5}}{5}}$ | $15\cdot a^2\sqrt{\dfrac{5+2\cdot\sqrt{5}}{5}}$ | $\dfrac{a}{2}\sqrt{\dfrac{25+11\cdot\sqrt{5}}{10}}$ | $\dfrac{5\cdot a^3}{2}\cdot\sqrt{\dfrac{47+21\cdot\sqrt{5}}{10}}$ |

El volumen de un poliedro regular es la tercera parte del producto de su área por la apotema.

En general el área de la cara de un poliedro regular viene dada por $A = kr^2$ donde $r$ es la apotema y $k$ es un número real. De manera similar se tiene para el número de diagonales de un polígono regular.

### *La función cuadrática y la geometría analítica*

En realidad la función cuadrática aparece el Teorema de Pitágoras para la distancia entre puntos. Recordemos que en todo triángulo rectángulo de catetos de medidas $a, b$ y de hipotenusa de medida $h$ se cumple que $h^2 = a^2 + b^2$ que es una función cuadrática en dos variable. Este resultado se generaliza a distancia entre puntos del espacio n-dimensional. De este hecho se obtiene la longitud de una curva, su curvatura etc. Generalizando las superficies cuadráticas son casos particulares de la forma cuadrática en n variables.

Uno de los resultados importantes es que la ecuación cuadrática permite demostrar la ortogonalidad entre la recta tangente en el punto de tangencia y una esfera en $R^n$.

## *La función cuadrática y teoría de números:*

### *La función cuadrática y los números perfectos*

Sabemos que

$$1 + 2 + \cdots + n = \frac{n(n+1)}{2}$$

Es decir, la función cuadrática $f(x) = \frac{x^2+x}{2}$ contiene la suma de los primeros n números naturales. Lo interesante es que si $a = 2^{p-1}(2^p - 1)$ es un número perfecto, vemos que

$$1 + 2 + \cdots + 2^p - 1 = \frac{2^p(2^p - 1)}{2} = 2^{p-1}(2^p - 1)$$

Nótese que en este caso $x = 2^{p-1} - 1$. De otra parte si hacemos $x = 2x + 1$ en la función $f$ entonces se obtiene

$$g(x) = 2x^2 + 3x + 1.$$

Esta función cuadrática es especial porque envía los números impares en números pares y los enteros pares en enteros impares. De hecho si $a$ es un número perfecto par, existe $n$ de tal manera que $2n^2 + 3n + 1 = a$.

**Definición.** A la parábola $f(x) = 2x^2 + 3x + 1$ la llamamos la parábola perfecta o de Elías.

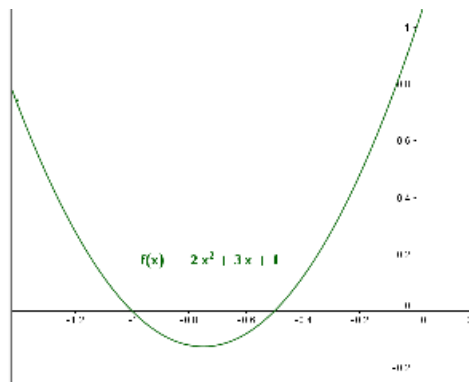

Si en $f$ hacemos $x = 2x$ obtenemos la función

$$h(x) = 2x^2 + x$$

Esta función cuadrática envía números pares en números pares e impares en impares. Y curiosamente no existe un entero $n$ de tal manera que $g(n)$ sea un número perfecto.

La suma de los primeros $n$ números impares (pares) viene dada por la función cuadrática:

$$f(n) = n^2 \; (f(n) = n^2 + n).$$

En general la suma

$$b + (b+d) + (b+2d) + \cdots + (b + (n-1)d) = \frac{dn^2 + (2b-d)n}{2} = f(n)$$

Nótese que en este caso

$$f(n) = \frac{(2dn + 2b - d)^2 - (2b-d)^2}{8d}$$

Si $d = b = 1$ se tiene

$$f(n) = \frac{(2n+1)^2 - 1}{8}$$

Esto explica porque el cuadrado de un número impar no puede ser un número perfecto.

**Estudiemos un poco la parábola de Elías.**

En primer lugar los ceros de $f(x) = 2x^2 + 3x + 1$ son -1 y -1/2. Entonces,

- El vértice se encuentra en $x = -\frac{3}{4}$, $f\left(-\frac{3}{4}\right) = -0.125 = -\frac{1}{8}$.
- Área encerrada por la parábola de Elías y el eje x.
$$A = 1/24$$
- Área encerrada por la parábola de Elías y el eje x entre -0.5 y 0
$$A = \frac{5}{24}$$
Nótese que la suma de las dos áreas calculadas es $\frac{1}{4}$.
- $f(\pi) = 2\pi^2 + 3\pi + 1 = 30.1639.. \equiv 30$
- $f(e) = 2e^2 + 3e + 1 = 23.9329 \ldots \equiv 24$
- $f(\varphi) = 2\varphi^2 + 3\varphi + 1 = 11.090 \equiv 11$
- DE otro lado, si tenemos
$$f(a) = 2a^2 + 3a + 1, f(b) = 2b^2 + 3b + 1 \; entonces,$$

$$f(a) - f(b) = (a-b)(2a + 2b + 3)$$

Calculando por separado vemos que:

- $\pi - e = 0.423310825130748$, $2(\pi + e) = 11.719748964097677$
- $\pi - \varphi = 1.52355866483989838846$, $2(\pi + \varphi) = 9.51925328467937617692$
- $e - \varphi = 1.10024783970915038536$, $2(e + \varphi) = 6.67263163441788017072$

Lo más interesante es:

$$2\pi + e = 9.0014671356386317 \equiv 9$$

$$2e + \pi = 8.5781563105078837 \equiv 8.6$$

$$2\pi + \varphi = 7.90121929592948132693 \equiv 8$$

$$2\varphi + \pi = 6.3776606310895829385 \equiv 6$$

$$2e + \varphi = 7.05459764566798532072 \equiv 7$$

$$2\varphi + e = 5.9543498059588349354 \equiv 6$$

$$2\pi + 3e = 14.4380307925567222 \equiv 14.5$$

$$2e + 3\pi = 14.8613416176874702 \equiv 15$$

$$2\varphi + 3\pi = 12.6608459382691694154 \equiv 13$$

$$3\varphi + 2\pi = 11.13728727342927102693 \equiv 11$$

De otro lado sabemos que

$$2x^2 + 3x + 1 = (x+1)(2x+1)$$

Los números terminados en 1, 3,5 son los que nos sirven para hallar números perfectos. Recordemos que si $ax^2 + bx + c = 0$ entonces por la formula cuadrática,

$$x = \frac{-b \pm \sqrt{b^2 - 4ac}}{2a} = \frac{-b \pm h}{2a}$$

Sea $x_1 = \frac{-b+h}{2a}$ entonces, $h = 2ax_1 + b = f'(x_1)$. Es decir, la raíz cuadrada del discriminante de la ecuación, es la derivada de la función cuadrática evaluada en el cero dado. Es decir, $\Delta = \left(f'(x_1)\right)^2$. Además nótese que el segundo cero es $x_2 = -\frac{ax_1+b}{2a}$.

En especial si hacemos $f(x) = 2x^2 + 3x + 1 - P = 0, f(x_1) = P$ en la parábola de Elías, encontramos que $x_2 = -\frac{3+2x_1}{2}$. Aplicando la fórmula cuadrática vemos que

$$x_1 = \frac{-3 + \sqrt{1+8P}}{4} \to P = \frac{(4x_1+3)^2 - 1}{8} = \frac{(f'(x_1))^2 - 1}{8}$$

Una consecuencia inmediata es que los números pares generados por esta función (en particular los perfectos) son divisibles al menos por $2^4$.

De otro lado algunos valores son:

| $X_1$ | $X_2 = -\frac{3+2x_1}{2}$ | P |
|---|---|---|
| 1 | -5/2 | 6 |
| 3 | -9/2 | 28 |
| 5 | -13/2 | 66 |
| 7 | -17/2 | 120 |
| 15 | -33/2 | 496 |
| 63 | -129/2 | 8128 |
| 4095 | -8193/2 | 33 550 336 |
| 1023 | -2046/2 | 2096128 |

De otro lado si hacemos,

- $x_1 = 2^l + 1 \to P = 2^l[2^{l+1} + 7] + 6$
- $x_1 = 2^l - 1 \to P = 2^l(2^{l+1} - 1)$
- $x_1 = 2n + 1 \to P = 8n^2 + 14n + 6$ (*)

Es decir, todo número perfecto par está en la parábola (*). Nótese además que si $P = 2^{p-1}(2^p - 1)$ tenemos que:

$$8P + 1 = (2^{p+1} - 1)^2$$

$$P = \frac{(2^{p+1}-1)^2 - 1}{8}$$

Por lo que Por otro lado si recordamos que

$$\sum_{i=1}^{n} i^2 = \frac{n(n+1)(2n+1)}{6}$$

Tenemos que

$$\sum_{i=0}^{n} i^2 = \frac{nf(n)}{6}, f(n) = (n+1)(2n+1)$$

Y por tanto,

$$f(n) = \frac{6}{n}\sum_{i=1}^{n} i^2 = \frac{6n^2}{n^3}\sum_{i=1}^{n} i^2 = 6n^2 \left(\frac{1}{n^3}\sum_{i=1}^{n} i^2\right)$$

Si $n$ es suficientemente grande $\frac{1}{n^3}\sum_{i=1}^{n} i^2 \to \frac{1}{3}$ por tanto

$$f(n) = 2n^2$$

Es decir, los números perfectos pares grandes son de la forma dada. Nótese que si $P$ es un número perfecto, $n = \left[\sqrt{\frac{P}{2}}\right] = \left[\frac{\sqrt{(2^{p+1}-1)^2-1}}{4}\right]$

### *Área en la parábola de Elías.*

Sabemos que $f(x) = 2x^2 + 3x + 1$, tenemos, $f(a) = 2a^2 + 3a + 1, f(b) = 2b^2 + 3b + 1$. La pendiente de la recta que pasa por los puntos $(a, f(a)), (b, f(b))$ viene dada por

$$m = 2a + 2b + 3$$

Y ecuación $y = (2a + 2b + 3)x + 1 - 2ab$. Calculemos,

$$\int_a^b y\,dx = \frac{b-a}{2}(f(a) + f(b))$$

Nótese que la figura formada es un trapecio y la expresión que aparece es el área de dicho trapecio.

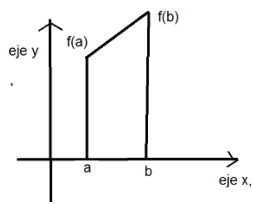

Ahora,

$$\int_a^b (2x^2 + 3x + 1)dx = \frac{b-a}{6}(2f(a) + 2f(b) + 4ab + 3a + 3b + 2)$$

Restando las áreas encontramos que el área encerrada por la recta y la parábola viene dada por:

$$A = \frac{(b-a)^3}{3}$$

### *La función cuadrática y la Conjetura de Goldbach.*

Comencemos con una propiedad interesante que cumplen los números enteros impares. Recordemos que si $p$ es impar entonces $p$ es de la forma $4k + 1$ o $4v - 1$ pero no de ambas a la vez.

**Teorema.** Sean $p, q$ números impares entonces

$$\frac{p+q}{2} \text{ es par (impar)} \Leftrightarrow \frac{p-q}{2} \text{ es impar (par)}$$

**Demostración** En efecto consideremos las siguientes posibilidades:

1. $p = 4k + 1, q = 4l - 1$ entonces

$$\frac{p+q}{2} = \frac{4k+1+4l-1}{2} = 2(k+l) \Leftrightarrow \frac{p-q}{2} = \frac{4k+1-4l+1}{2} = 2(k-l) + 1$$

2. $p = 4k + 1, q = 4l + 1$ entonces

$$\frac{p+q}{2} = \frac{4k+1+4l+1}{2} = 2(k+l)+1 \Leftrightarrow \frac{p-q}{2} = \frac{4k+1-4l-1}{2} = 2(k-l)$$

3. $p = 4k-1, q = 4l-1$ entonces

$$\frac{p+q}{2} = \frac{4k-1+4l-1}{2} = 2(k+l)-1 \Leftrightarrow \frac{p-q}{2} = \frac{4k+1-4l+1}{2} = 2(k-l)$$

4. $p = 4k-1, q = 4l+1$

$$\frac{p+q}{2} = \frac{4k-1+4l+1}{2} = 2(k+l) \Leftrightarrow \frac{p-q}{2} = \frac{4k-1-4l-1}{2} = 2(k-l)-1$$

Así considerando todas las posibilidades se tiene el resultado.

**Un resultado simple**

Sean $p, q$ números impares, hagamos

$$M = \frac{p+q}{2}, \ I = \frac{p-q}{2}$$

Si suponemos $M$ par entonces $I$ es impar. Sea $M = 2n$ podemos escribir entonces

$$p + q = 2n, p - q = 2I$$

Sumando y restando encontramos p=$2n + I, q = 2n - I$. También tenemos que $p = 2I + q$. Ahora el teorema de Dirichlet afirma si $a, b$ son primos relativos entonces la sucesión $p = an + b$ contiene infinitos números primos, en particular si elegimos los números 2 con $q$ primo dado y $3 \leq q \leq 2n-1$. Entonces la sucesión

$$p = 2m + q$$

contiene infinitos números primos. Elijamos $m$ un impar $I$ para el cual $p$ es primo, es decir,

$$p = 2I + q$$

Devolviendo los pasos vemos que

$$p - I = q + I = 2l$$

ya que la suma y resta de impares es par. Vemos también que

$$p = 2l + I, q = 2l - I, p + q = 4l$$

De estas igualdades y las anteriores y notando que $q$ es de la forma $2n - I$ vemos que se puede escoger $l = n$ y el resultado se tiene. Un resultado similar se obtiene si $M$ es impar. Se tiene entonces el siguiente resultado:

**Teorema** Para cada $n \geq 2$ existe un impar $I < n$ de tal manera que

$$p = 2n + I, q = 2n - I$$

son números primos

Como consecuencia se tiene el siguiente resultado.

**Goldbach-Elías** Todo par $n \geq 4$ se puede escribir como la suma de dos números primos.

### Parábolas

Sean $p > q$ números impares y consideremos la parábola

$$f(x) = x^2 - (p+q)x + pq$$

Notemos que los cortes con el eje $x$ son $x = p, x = q$. Además, el vértice ocurre en

$$x = -\frac{b}{2a} = \frac{p+q}{2} = M$$

y para tal caso

$$f(\tfrac{p+q}{2}) = -\frac{(p-q)^2}{4} = -I^2 \quad , I = \frac{p-q}{2}$$

Si $M$ es par (impar) $I$ es impar (par).

## Algunos cálculos

Calculemos el área $A_s$ encerrada por el eje $x$ y la parábola entre los puntos $x = p, x = q$ viene dada por

$$A_s = \int_q^p (-x^2 + (p+q)x - pq)\, dx = \frac{(p-q)^3}{6} = \frac{4}{3} I^3$$

El área del rectángulo $A_r$ de lados $p - q, I^2$ viene dada por

$$A_r = (p-q)I^2 = 2I^3$$

El área del triángulo $A_t$ de base $p - q$ y altura $I^2$ viene dada por

$$A_t = \frac{(p-q)I^2}{2} = \frac{(p-q)^3}{8} = I^3$$

Area bajo la curva

$$A_s = \int_0^q (x^2 - (p+q)x + pq)\, dx = \frac{1}{6} q^2 (3p - q)$$

Tenemos las siguientes relaciones

$$\frac{A_r}{A_s} = \frac{3}{2} \Leftrightarrow A_r = \frac{3}{2} A_s$$

$$\frac{A_r}{A_t} = 2 \Leftrightarrow A_r = 2A_t$$

$$\frac{A_s}{A_t} = \frac{4}{3} \Leftrightarrow A_s = \frac{4}{3} A_t$$

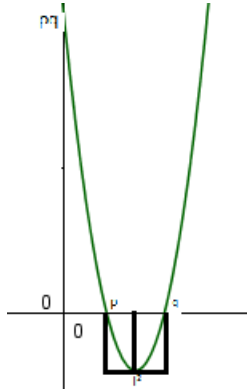

Nótese que:

$$p = 2n + I, \ q = 2n - I$$

**Triángulos rectángulos y números primos**

De otra parte note que en el triángulo rectángulo de lados $2n$ e $I$ con números primos

$$p = 2n + I, q = 2n.I \Leftrightarrow p + q = 4n, p - q = 2I$$

Elevando al cuadrado tenemos

$$(p + q)^2 = p^2 + 2pq + q^2$$
$$(p - q)^2 = p^2 - 2pq + q^2$$

Sumando

$$(p + q)^2 + (p - q)^2 = 16n^2 + 4I^2 \Leftrightarrow \frac{(p + q)^2 + (p - q)^2}{4} = 4n^2 + I^2$$

Nótese que como $((2n, I)) = 1$ entonces $((4n^2, I^2)) = 1$ por lo que el número

$$H^2 = (2n)^2 + I^2$$

Puede ser un número primo o un primo al cuadrado.

Ejemplo: 12+5=17, 12-5=7. Aquí n=6, luego 4(36)+25=169=13$^2$ También, 12+7=19, 12-7=5. Por lo que 4(36)+49= 144+49=193 que es primo.

Sean $p = 2n + I, q = 2n - I$ números primos. Entonces,

$$h = 4n^2 + I^2$$

puede ser un número primo.

$$H = (2n)^{2l} + I^{2l}, l \geq 1$$

puede ser un número primo.

Nótese que

$$(p+q)^{2l} + (p-q)^{2l} = (4n)^{2l} + (2I)^{2l} \Leftrightarrow \frac{(p+q)^{2l} + (p-q)^{2l}}{2^{2l}} = (2n)^{2l} + I^{2i}$$

## *La función cuadrática y los números metálicos.*

Recordemos que los números metálicos fueron descubiertos por la matemática argentina Vera W. de Spinadel allá cerca del año 2000. Estos resultan de resolver una ecuación cuadrática. Consideremos los números

$$x_1 = \frac{1+\sqrt{m}}{2}, \quad x_2 = \frac{1-\sqrt{m}}{2}, \ m \in Z^+$$

Que conducen a la ecuación cuadrática

$$x^2 - x - \frac{m-1}{4} = 0 \ (t)$$

Ahora, la expresión $\frac{m-1}{4}$ es un entero si y solo si $m = 4n + 1, n \in Z^+$, es decir, la ecuación (t) se reduce a

$$x^2 - x - n = 0 \quad (t')$$

con $n$ entero positivo. De igual forma los números

$$x_1 = \frac{1+\sqrt{m}\,i}{2}, \quad x_2 = \frac{1-\sqrt{m}\,i}{2}, \ m \in Z^+$$

conducen a la ecuación

$$x^2 - x + \frac{1+m}{4} = 0 \ (f)$$

De manera similar $\frac{1+m}{4}$ es un entero si y solo si $m = 4n - 1$. Y la ecuación se reduce a

$$x^2 - x + n = 0 \quad (f')$$

Con $n$ entero positivo. Lo interesante de la discusión anterior muestra que las ecuación (t), (f) tienen soluciones en donde el radical es un número impar de la forma $4n + 1$ o de la forma $4n - 1$ según, el caso. Es decir, la mitad de los números impares producen soluciones complejas y la otra mitad soluciones reales. Como sabemos un número primo es de esta forma. Nótese en la ecuación (t) se obtienen soluciones reales y distintas, si el número $m$ no es cuadrado perfecto, la soluciones están en el cuerpo cuadrático $Q(\sqrt{m})$. Cuando $m$ es un número primo impar, se observa una propiedad interesante: El número primo caracteriza la ecuación y aparece en la solución.

Nótese en la tabla de números metálicos:

| p | q | Ecuación | Solución positiva | Número metálico | |
|---|---|----------|-------------------|-----------------|---|
| 1 | 1 | $x^2-x-1=0$ | $\frac{1+\sqrt{5}}{2}$ | $\sigma_{1,1} = \frac{1+\sqrt{5}}{2}$ | (número de oro) |
| 2 | 1 | $x^2-2x-1=0$ | $\frac{2+2\sqrt{2}}{2}$ | $\sigma_{2,1} = 1+\sqrt{2}$ | (número de plata) |
| 3 | 1 | $x^2-3x-1=0$ | $\frac{3+\sqrt{13}}{2}$ | $\sigma_{3,1} = \frac{3+\sqrt{13}}{2}$ | (número de bronce) |
| 4 | 1 | $x^2-4x-1=0$ | $\frac{4+2\sqrt{5}}{2}$ | $\sigma_{4,1} = 2+\sqrt{5}$ | |

Los números impares no son solamente primos sino que son de la forma $4n + 1$, excepto por el número de la plata donde aparece el número dos que también es primo. Esto no quiere decir que todo el número impar bajo la cantidad subradical tiene que ser primo. De otro lado los números metálicos deber ser ceros de la ecuación

$$x^2 - x - n = 0 \quad (t')$$

y no necesariamente de la ecuación

$$x^2 - nx - 1 = 0 \quad (t')$$

Como se muestra en la tabla.

De otro lado se sabe que el número de oro es el más irracional de los números metálicos, sin embargo vemos que:

$$1.6 < \varphi < 1.7$$

Por tanto es posible que en las medidas donde aparece el número áureo, la aproximación esté entre estos dos valores.

Supongamos que los números
$$x_1 = \frac{1+\sqrt{m}}{2}, \quad x_2 = \frac{1-\sqrt{m}}{2}$$

son números enteros, entonces $m$ tiene que ser un cuadrado perfecto, en particular, como $m$ es impar su raíz también, digamos, $\sqrt{m} = 4k+1$ por lo que

$$2x_1 = 4k+2 \leftrightarrow x_1 = 2k+1 \leftrightarrow x_2 = -2k$$

Que son las raíces de la ecuación

$$x^2 - x - 2k(2k+1) = 0$$

Vemos que $n = 2k(2k+1)$ por lo que $m = 4n+1 = (4k+1)^2$

## La función cuadrática y los p-números.

En esta parte consideramos los números de la forma: $pppppppppppp$ ($t$ veces), $1 \leq p \leq 9$

Los cuales llamaremos p-números. Como veremos esos números son de especial importancia. Escribiremos $p \times t$ para indicar el p-número. En particular si $p$ es un dígito escribimos $p \times 1$. Es decir,

$$pppppppppppp \equiv p \times t$$

Si $p$ es un entero de más de dos cifras que no es un p-número le hacemos corresponder el p-número $a_1 \times 1$ donde $a_1 \neq 0$ es el digito de las unidades. Si $a_1 = 0$ asociamos al número entero el pnúmero $a_1 + 1 \times 1$.

También podríamos asociar a un entero cualquiera el digito de sus unidades. Esto obviamente define una relación de equivalencia en el conjunto de los números enteros. También podríamos definir en el conjunto de los números enteros la relación de sumar sus dígitos hasta que se obtenga un digito, por ejemplo, si un entero la suma de sus dígitos da 124 le asociamos el número 7, si por ejemplo la suma de sus dígitos da 80 le asociamos el cero, o si la suma de sus dígitos es 12349998 sumamos una vez más da 45, sumando da 9 y este el número que le asociamos.

A todo p-número le podemos asociar una parábola, si $p \times q$ podemos formar la parábola

$$x^2 - (p+q)x + pq$$

La cual tiene su conjugada

$$x^2 + (p+q)x + pq$$

En conclusión a cada entero le podemos asociar un p-número y a cada p-número una parábola.

## La función cuadrática y la física.

La función cuadrática tiene muchas aplicaciones en la física, quizá el ejemplo más simple es el lanzamiento de proyectiles, sabemos que si un objeto es lanzado con una ángulo de elevación β con una velocidad inicial $V_0$ entonces, el movimiento en el eje $y$ viene dado por la fórmula $y = -\frac{gt^2}{2} + V_0 sen(\beta)t$. Además, como el movimiento en el eje horizontal es uniforme, $x = V_0 \cos(\beta)t$ eliminando el parámetro encontramos que

$$y = -\frac{gx^2}{2V_0^2 \cos^2(\beta)} + \tan(\beta)x$$

Que es una función cuadrática.

De hecho muchas expresiones tales como

$$E = mc^2, K = \frac{1}{2}mv^2$$

Son parábolas. Es común también encontrar expresiones muy importantes en física en donde se divide por el cuadrado de algo. Por ejemplo, la ley de atracción de Newton para dos masas separadas por una distancia $d$ viene dada por: $F = \frac{GMm}{d^2}$ donde $G$ es la constante universal por el producto de las masas y $d$ es la distancia al cuadrado de sus centros.

La fuerza centrífuga viene dada por $F = \frac{mv^2}{d^2} = ml^2, l = \frac{v}{d}$

## PARA FINALIZAR: ALGUNAS OBSERVACIONES SOBRE EL NÚMERO $\varphi$

Sabemos que el número $\varphi$ es un cero de la ecuación cuadrática $x^2 - x - 1 = 0$.

Resumamos algunas de las propiedades vistas hasta ahora de este maravilloso número.

Recordemos que
$$\varphi = \frac{1+\sqrt{5}}{2},$$

Luego,

Sabemos que $\varphi = \dfrac{1+\sqrt{5}}{2}, \overline{\varphi} = \dfrac{1+\sqrt{5}}{2}$, entonces, 2

1. $\varphi + \overline{\varphi} = 1$
2. $\varphi \cdot \overline{\varphi} = -1$
3. $\varphi^2 + \overline{\varphi}^2 = 3$
4. $\varphi^2 = \varphi + 1$
5. $\varphi^2 + \varphi\overline{\varphi} = \varphi$
6. $\varphi^3 = 2\varphi + 1$
7. $\varphi^4 = 3\varphi + 2$
8. $\varphi^5 = 5\varphi + 3$
9. $\varphi^6 = 8\varphi + 3$

Notemos que en general, si $\varphi^n = a\varphi + b$ entonces,

1. $\overline{\varphi}^n = a\overline{\varphi} + b$ por lo que $\varphi^n + \overline{\varphi}^n = a + 2b$ y $\varphi^n - \overline{\varphi}^n = -\sqrt{5}a$
2. $\varphi^{n+1} = (a+b)\varphi + a$

O también, si $\varphi^{n-1} = A\varphi + B$ tenemos que $\varphi^n = (A+B)\varphi + A$ y $\overline{\varphi}^n = (A+B)\overline{\varphi} + A$. Por lo que
$$\varphi^n + \overline{\varphi}^n = 3A + B$$
$$\varphi^n - \overline{\varphi}^n = \sqrt{5}(A + B)$$

Si hacemos, $A = F_{n-1}, B = F_{n-2}$ números de Fibonnaci encontramos,
$$\varphi^n + \overline{\varphi}^n = 3F_{n-1} + F_{n-2} = F_{n+1} + F_{n-1}$$
$$\varphi^n - \overline{\varphi}^n = \sqrt{5}(F_{n-1} + F_{n-2}) = \sqrt{5}F_n$$